\documentclass[a4paper, oneside, 10pt]{article}
\usepackage[english]{babel}
\usepackage{amsfonts}
\usepackage{amssymb}
\usepackage{graphics}
\usepackage{latexsym}
\begin{document}

\title{{\bf {\Large{A note on Bayesian nonparametric priors derived from exponentially tilted Poisson-Kingman models.}}}\footnote{{\it AMS (2000) subject classification}. Primary: 60G57. Secondary: 60G09.}}
\author{\textsc {Annalisa Cerquetti}\footnote{Corresponding author.
E-mail: {\tt annalisa.cerquetti@unibocconi.it}}\\
  \it{\small Bocconi University, Milano, Italy }}
\newtheorem{teo}{Theorem}
\date{}
\maketitle{}

\begin{abstract}
We derive the class of normalized generalized Gamma processes from Poisson-Kingman models (Pitman, 2003) with tempered $\alpha$--stable mixing distribution. Relying on this construction it can be shown that in Bayesian nonparametrics, results on quantities of statistical interest under those priors, like the analogous of the Blackwell-MacQueen prediction rules or the distribution of the number of distinct elements observed in a sample, arise as immediate consequences of Pitman's results.\\

\noindent{\it Keywords}: Exchangeable random partitions;  Exponential tilting; Inverse Gaussian density; Random probability measures; Tempered stable laws.

\end{abstract}

\section{Introduction}
In Lijoi, Mena and Pr\"unster (2005) the normalized inverse Gaussian (N-IG) process has been introduced as an alternative to the Dirichlet process to be used in Bayesian nonparametric mixture modeling. 
By mimicking Ferguson's (1973) famous construction of the Dirichlet process, the authors define a random discrete probability measure $P$, on a Polish space $(S, \mathcal{S})$, whose finite dimensional distributions have the multivariate law of a vector of $n$ independent r.v.'s with inverse Gaussian distribution divided by their sum. Even if the authors observe that a N-IG prior, with non-atomic parameter measure, belongs to the class of {\it species sampling models} (Pitman, 1996), and more exactly
to {\it Poisson-Kingman models} (Pitman, 2003), they do not derive such process as an element of the previous families,
but obtain, independently of  Pitman's results, both the analogous of the Blackwell-MacQueen prediction rules and the distribution of the number of distinct values observed in a sample.

Here we show how the larger class of normalized generalized Gamma processes (N-GG), already considered in James (2002), may be derived from Poisson-Kingman models for random partitions of the positive integers. In particular these processes arise as random discrete probability measures whose ranked atoms follow a Poisson-Kingman distribution derived from an $\alpha-$stable law with mixing distribution the {\it exponentially tilted} version of the stable density. It follows that N-GG priors induce {\it exchangeable Gibbs  partitions  of type $\alpha$}, (see Gnedin and Pitman, 2005), hence, even if Pitman (2003) is not directly concerned with applications in Bayesian nonparametrics, distributional results on quantities of statistical interest  under those priors, arise as straightforward consequences of Pitman's results.

The paper is organized as follows. In Section 2 we recall the definition of exponentially tilted Poisson-Kingman models derived from a stable law of index $\alpha$.  In Section 3 we exploit results in Pitman (2003) to derive a general expression for corresponding exchangeable partition probability function (EPPF) and analogous of the Blackwell-MacQueen prediction rules. Finally, in Section 4, the distribution of the number of blocks and its asymptoyic behaviour are obtained from results in Pitman (2002, 2003) and Gnedin and Pitman (2005). 

\section{Exponentially tilted $\alpha-$stable Poisson-Kingman models}
It is well known that, given a law $Q$ on the space $\mathcal{P}_1^\downarrow$  of decreasing sequences of positive numbers with sum 1 and a law $H(\cdot)$ on a Polish space $(S, \mathcal {S})$, a random discrete probability measure $P$ on $\mathcal{S}$ may be defined by $P(\cdot)=\sum_{i=1}^\infty P_i \delta_{X_i}(\cdot)$, for $X_i$ iid $\sim H(\cdot)$ and $(P_i) \sim Q$. Generalizing Kingman's (1975) construction of the Dirichlet process as a Gamma process with independent increments divided by the sum, Jim Pitman $-$ in a stimulating paper available on his web page since 1995 and published in 2003 $-$ introduces a large class of random discrete probability measures deriving the law of the atoms, in decreasing order, from the ordered points of a homogeneous Poisson process on $(0, \infty)$ with given L\'evy density divided by their sum.\\\\
{\bf Definition 1.} [Pitman (2003; Def. 3)] { Let $P_i=(J_i/T)$ be a ranked discrete distribution derived from the ranked points of a Poisson process with L\'evy density $\rho$  of random lenghts $J_1 \geq J_2 \geq \cdots \geq 0$ by normalizing their lenghts by their sum which is $T$. 
The law $Q$ on $\mathcal{P}_1^{\downarrow}$ of the sequence $(P_i)$ will be called the {\it Poisson-Kingman distribution with L\'evy density $\rho$}, and denoted PK$(\rho)$. Denote by  $PK(\rho|t)$ the regular conditional distribution of $(P_i)$ given $(T=t)$ constructed above. For a probability distribution $\gamma$ on $(0,\infty)$, let
\begin{equation}
\label{PKMIX}
PK(\rho, \gamma):=\int_0^\infty PK(\rho|t) \gamma(dt)
\end{equation}
be the distribution on the space $\mathcal{P}_1^{\downarrow}$. Call $PK(\rho, \gamma)$ the {\it Poisson-Kingman distribution with L\'evy density $\rho$ and mixing distribution $\gamma$}}.\\\\
{\bf  Remark 1.} In James, Lijoi and Pr\"unster (2005) a very large class of Normalized Random Measure (NRMs) based on a more complex generalization of Kingman's construction, has been introduced and deeply studied in a Bayesian nonparametric perspective.  As the same authors point out, Pitman's $PK(\rho, \gamma)$ models provide an important extension of {\it homogenous} NMRs. Those model contain, among the others, the two-parameter Poisson--Dirichlet distribution, $PD(\alpha, \theta)$, for $0\leq \alpha < 1$ and $\theta > -\alpha$, which is the law of the ranked atoms of the well-known extension of the Dirichlet process introduced in Pitman and Yor (1997). Pitman (2003) shows this family corresponds to a 
Poisson-Kingman model derived from a stable law of index $\alpha$ with mixing distribution $\gamma(t)=\frac{\Gamma(\theta+1)}{\Gamma (\theta /\alpha +1)} t^{-\theta} f_{\alpha, \delta}(t)$, for $f_{\alpha,\delta}(t)$ the density of the stable law.  (See also Pitman, 2002, for an exaustive account of this family of distributions on $\mathcal{P}_1 ^{\downarrow}$.)\\\\ 
In what follows we shall actually focus on models $PK(\rho_\alpha, \gamma)$ where $\rho_\alpha$ is the L\'evy density of a stable density of index $\alpha \in (0,1)$. The reason lies in Theorem 8, Section 5.3, in Pitman (2003) and will be clarified later.\\\\
First recall that, given a strictly positive r.v. $T$, with  density $f_T$ and Laplace transform 
$
E(e^{-\lambda T})= e^{-\psi(\lambda)} =\int_{0}^{\infty} e^{-\lambda t}f_T (t)dt,
$
where, according to  L\'evy-Kintchine formula, $\psi(\lambda)=\int_{0}^{\infty} (1 - e^{-\lambda s})\rho(ds)$ is the Laplace exponent, $f_T$  is uniquely determined by its unique L\'evy density $\rho$.  Now, one of the basic operations that leads to the larger class of $PK(\rho, \gamma)$ models introduced in Definition 1 is given by {\it exponential tilting} (see Pitman, 2003, Sec. 4.2). Tilting a positive random variable $T$ is performed by multiplying its density by a  factor $exp\{\psi(\lambda) -\lambda t\}$, or its L\'evy density by a factor $exp\{-\lambda t\}$. For our purposes it is worth to recall the definition of {\it tempered stable law}, first introduced in Tweedie (1984), also called generalized Gamma distributions in Brix (1999), (see Barndorff-Nielsen and Shephard, 2001, for a comprehensive account).\\\\
{\bf Definition 2.} [Tempered $\alpha$-stable law] Let $f_{\alpha, \delta}(t)$, $\alpha \in (0,1)$, $\delta \in (0,\infty)$, denote the probability density function of a positive $\alpha$-stable law with Laplace exponent $\psi_\alpha(\lambda)=\delta(2\lambda)^{\alpha}$. Apart from $\alpha=\frac 12$ and $\alpha=\frac 13$, for a general $\alpha \in (0,1)$ explicit expressions of this density are known only in the form of series representations:
$$
f_{\alpha, \delta}(t)=\frac{1}{2\pi} \delta^{-1/\alpha} \sum_{\xi=1}^{\infty} (-1)^{\xi-1}\sin(\xi \pi \alpha)\frac{\Gamma(\xi \alpha +1)}{\xi!} 2^{\xi \alpha+1}(t/\delta^{1/\alpha})^{-\xi \alpha-1}).
$$
By the change of variable $\lambda=\frac{\gamma^{\frac 1\alpha}}{2}$, $\gamma \in [0,\infty)$, exponential tilting $f_{\alpha, \delta}(t)$ with $exp\{\delta\gamma - \frac \gamma 2^{\frac 1\alpha}t\}$, gives the density of a {\it tempered stable law} of parameters $(\alpha, \delta, \gamma)$, 
\begin{eqnarray}
&&f_{\alpha, \delta, \gamma}(t)= e^{\delta \gamma -\frac 12 \gamma^{1/\alpha}t} f_{\alpha,\delta}(t)\nonumber.
\end{eqnarray}
Although this density has no explicit expression, corresponding Laplace exponent and L\'evy density are known to be as follows:
\begin{equation}
\label{LPTSA}
{\psi_\alpha^e}(\lambda)=-\delta\gamma +\delta(\gamma^{\frac 1 \alpha} +2\lambda)^\alpha
\quad {\mbox {and}} \quad
{\rho}_\alpha^e(s)=\delta 2^\alpha \frac{\alpha}{\Gamma(1-\alpha)}s^{-1-\alpha}e^{-\frac 12\gamma^{\frac1 \alpha}s}.
\end{equation}
{\bf Example 1.} [Inverse Gaussian law] The class of tempered $\alpha$-stable laws contains the inverse Gaussian law. In fact, for $\alpha=\frac 12$ the stable density has the following  explicit form 
\begin{eqnarray}
\label{stable12d}
&&f_{1/2, \delta}(t)=\frac{\delta}{\sqrt{2\pi}}t^{-\frac 32} e^{-\frac {\delta^2}{2t}}\nonumber, 
\end{eqnarray}
and corresponding Laplace exponent $\psi(\lambda)=\delta\sqrt{2\lambda}$. By exponential tilting with $\lambda=\frac {\gamma^2}{2}$, the density of a tempered $\frac 12-$stable  law results
\begin{equation}
\label{igdens}
f_{\delta, \gamma}(t)= 
\frac{\delta}{\sqrt{2\pi}} e^{\delta \gamma}{t^{-\frac 32}} 
\exp\left\{-\frac 12 \left(\delta^2t^{-1}+\gamma^2 t\right)\right\},
\end{equation}
for $\delta \in (0,\infty)$ and $\gamma \in [0,\infty)$, which is well-known to be the density of an inverse Gaussian $(\delta, \gamma)$ law (see e.g. Seshadri, 1993). Corresponding Laplace exponent and  L\'evy density easily follow from (\ref{LPTSA}):
$$
\psi_{\frac 12}^e(\lambda)=-\delta\gamma +\delta(\gamma^2 +2\lambda)^{\frac 12}
\qquad\mbox{and}\qquad
\rho_{\frac 12}^e(s)=\frac {\delta}{\sqrt{2\pi}}s^{-\frac 32}e^{-\frac s2 \gamma^2}.
$$\\\\
It is well known that normalized generalized Gamma priors, as introduced in Bayesian nonparametric context, (see James, 2002) corresponds to random discrete probability measures $P(\cdot)=\sum_{i=1}^\infty P_i \delta_{X_j}(\cdot)$ whose ranked atoms $(P_i)$ follow a $PK(\rho_\alpha^e)$ distribution. Nevertheless in Section 4.2 (cfr. eq. (46)) Pitman shows that if $\rho^e$ is the tilted version of the L\'evy density of $T$, a model $PK(\rho^e)$ is equivalent to a model $PK(\rho, \gamma^e)$ where $\rho$ is the L\'evy density of $T$ and the mixing distribution, $\gamma^e$, is the tilted version of the density of $T$. This implies, for example, that  the normalized inverse Gaussian prior of Lijoi, Mena and Pr\"unster (2005) corresponds to a random discrete distribution whose ranked atoms have $PK(\rho_\frac 12, \gamma_{\frac 12}^e)$ distribution.  Relying on the previous considerations we are now in a position to introduce the following definition:\\\\
{\bf Definition 3.} Let $PK(\rho_\alpha)$ be a Poisson-Kingman model derived from an $\alpha$-stable law, and let $\gamma_\alpha^e(t)$, for $\alpha \in (0,1)$ denote the density of a tempered $\alpha-$stable law. We call the family of distributions $PK(\rho_\alpha, \gamma_\alpha^e) \equiv PK(\rho_\alpha^e)$ on $\mathcal{P}_1^{\downarrow}$ {\it exponentially tilted $\alpha-$stable Poisson-Kingman models}. 

\section{EPPF and predictive distributions}
In this section we derive distributional results for quantities of statistical interest in Bayesian nonparametric modeling under normalized generalized Gamma priors,  exploiting Pitman's results for Poisson-Kingman models. First recall that, from Kingman's theory of exchangeable random partitions (Kingman, 1978), sampling from a random discrete distribution $P$, induces a random partition $\Pi$ of the positive integers $\mathbb{N}$, by the exchangeable equivalence relation
$i \approx j \Leftrightarrow X_i=X_j$, that is to say two positive integers $i$ and $j$ belong to the same block of $\Pi$ if and only if $X_i=X_j$, where $X_i|P$ are iid $\sim P$. It follows that, for each restriction $\Pi_n=\{A_1,\dots, A_k\}$ of $\Pi$ to $[n]=\{1,\dots, n\}$, and for each $n=1,2,\dots$,  $Pr(\Pi_n=\{A_1,\dots, A_k\})=p(n_1,\dots, n_k)$, where, for $j=1,2,\dots, k$, $n_j=|A_j|\geq 1$ and $\sum_{j=1}^k n_j=n$, for some non-negative symmetric function $p$ of finite sequences of positive integers called the {\it exchangeable partition probability function} (EPPF) determined by $\Pi$.  
In Hansen and Pitman (2000) it is shown that an infinite exchangeable sequence $(X_n)$ admits prediction rules of the form 
\begin{equation}
\label{ssm}
Pr(X_{n+1}\in \cdot|X_1,\dots, X_n)=\sum_{j=1}^{K_n} p_{j,n} \delta_{X^*_j}(\cdot) +q_nH(\cdot)
\end{equation}
where  $X_j^*$, for $1\leq j\leq K_n$ are distinct values, in order of appearance, in $(X_1, \dots, X_n)$, $p_{j,n}$ and $q_n$ are non-negative product measurable functions of $(X_1,\dots, X_n)$, and $H(\cdot)$ is a non atomic probability measure on $\mathcal{S}$, if and only if $p_{j,n}=p({\bf n}^{j+})/p({\bf n})$ and $q_n=p({\bf n}^{k+1})/p({\bf n})$, where $p({\bf n}):=p(n_1,\dots, n_k)$, $p({\bf n}^{j+}):=p(n_1,\dots,n_j+1,\dots, n_k)$, and $p({\bf n}^{k+1}):=p(n_1,\dots,n_k, 1)$. Exchangeable sequences admitting predictive distributions of this form are termed {\it species sampling sequences}, and their directing measures $P$ are called {\it species sampling models} (Pitman, 1996).\\\\
By construction,  random discrete distributions derived by $PK(\rho, \gamma)$ models belong  to this class, therefore, to obtain from (\ref {ssm}) general expressions for the predictive distributions, one just need to know the EPPFs. Indeed  Pitman (2003) provides a thorough characterization of the laws $PK(\rho)$ on $\mathcal{P}_1^{\downarrow}$ via their corresponding EPPFs. Specifically, according to Corollary 6, Sec. 3, for some random partition  $\Pi_n=\{A_1,\dots, A_k\}$ of $[n]=\{1,\dots, n\}$, with block sizes $|A_i|=n_i$ for $i=1,\dots, k\leq n$, the EPPF associated with each $PK(\rho)$ is given by 
\begin{equation}
\label{EPPF}
p_K(n_1,\dots, n_k):=\frac{(-1)^{n-k}}{\Gamma(n)}\int_0^\infty \lambda^{n-1} e^{-\psi(\lambda)d\lambda}\prod_{i=1}^k \psi_{n_i}(\lambda) d\lambda,
\end{equation}
where $\psi(\lambda)$ is the Laplace exponent determined by $\rho(\cdot)$ and, for $m=1,\dots,n$, $\psi_m(\lambda):=\frac{d^m}{d\lambda^m}\psi(\lambda)$. It follows that, having at hand the equivalence stated in Definition 3., the EPPF of a $PK (\rho_\alpha, \gamma_\alpha^e)$ model can be easily deduced from (\ref{EPPF}) by substituting in  $\psi(\lambda)$ the Laplace exponent of the tempered $\alpha-$stable law given in (\ref {LPTSA}).\\\\
{\bf Proposition 1.} {\it Let $\gamma_\alpha^e(t)$ denote the density of a tempered   $\alpha$-stable law, for $\alpha \in (0,1)$, 
then the exchangeable partition probability function induced by an exchangeable sequence $(X_n)$ whose directing measure  $P$ has ranked  atoms following a $PK(\rho_\alpha, \gamma_\alpha^e)$ distribution, results  
\begin{equation}
\label{EPPFPKA}
p(n_1,\dots,n_k)=\frac{e^{\delta\gamma} \delta^k \alpha^k 2^n }{\Gamma(n)} \prod_{j=1}^k (1-\alpha)_{n_j -1 \uparrow}\int_{0}^{\infty} \lambda^{n-1} \frac{e^{-\delta (\gamma^{\frac 1\alpha}+2\lambda)^\alpha}}{(\gamma^{\frac 1\alpha}+2\lambda)^{n-k\alpha}}d\lambda,
\end{equation}
where  $n_j-{1\uparrow}$ stands for the usual notation of rising factorials $(x)_{n \uparrow}=x(x+1)(x+2)\cdots(x+n-1)$}.\\\\
{\it Proof}:\quad By equation (\ref{LPTSA}) the Laplace exponent of $\gamma_{\alpha}^e$ is given by  $\psi_\alpha^e(\lambda)=-\delta\gamma +\delta(\gamma^{\frac 1\alpha}+2\lambda)^{\alpha}$,
hence
\begin{eqnarray}
\psi_m(\lambda):=\frac{d^m}{d\lambda^m}\psi(\lambda)=\delta 2^{m}(\gamma^{\frac 1\alpha}+2\lambda)^{\alpha-m} (-1)^{m-1}\alpha \prod_{i=1}^{m-1}(\alpha -i)\nonumber.
\end{eqnarray}
By substitution in (\ref {EPPF})
\begin{eqnarray}
&&p(n_1,\dots,n_k)=
\frac{(-1)^{n-k}}{\Gamma(n)} \int_0^\infty \lambda^{n-1} e^{\delta\gamma -\delta(\gamma^{\frac 1\alpha}+2\lambda)^{\alpha}}\prod_{j=1}^{k}\delta 2^{n_j}\frac{(-1)^{n_j-1}\alpha \prod_{i=1}^{n_j-1}(\alpha -i)}{(\gamma^{\frac 1\alpha}+2\lambda)^{n_j-\alpha}}\nonumber
\end{eqnarray}
which reduces to 
\begin{eqnarray}
&&\frac{\delta^k\alpha^k 2^n e^{\delta\gamma}}{\Gamma(n)}   \prod_{j=1}^{k}( 1-\alpha)_{n_j-1\uparrow} \int_0^\infty \lambda^{n-1} \frac{e^{-\delta(\gamma^{\frac 1\alpha}+2\lambda)^{\alpha}}}{\prod_{j=1}^k (\gamma^{\frac 1\alpha}+ 2\lambda)^{n_j-\alpha}}d\lambda, \nonumber
\end{eqnarray}
and the result follows. \hspace{11.5cm}$\square$\\\\
Relying on the previous result, the general expression for predictive distributions induced by exponentially tilted $\alpha$--stable Poisson-Kingman models, easily follows.\\\\
{\bf Corollary 1.} {\it  An exchangeable sequence $(X_n)$ whose directing measure $P$ has ranked atoms $(P_i)$ with distribution $ PK(\rho_\alpha, \gamma_\alpha^e)$, has predictive distributions of the form $(4)$ for
\begin{eqnarray}
\label{weights}
&&p_{j,n}({\bf n})= \frac{2}{n} \frac{\eta_{n+1, k}}{ \eta_{n,k}}(n_j -\alpha)
\qquad {\mbox{and}}\qquad
q_n({\bf n})=\frac{2}{n} \frac{\eta_{n+1, k+1}}{ \eta_{n,k}}\alpha \delta
\end{eqnarray}
where}
\begin{eqnarray}
\label{vnk}
&&\eta_{n, k}=\int_{0}^{\infty} \lambda^{n-1} \frac{e^{-\delta (\gamma^{\frac 1\alpha}+2\lambda)^\alpha}}{(\gamma^{\frac 1\alpha}+2\lambda)^{n-k\alpha}}d\lambda.
\end{eqnarray}\\\\
{\bf Example 2.} [Normalized Inverse Gaussian process] Specializing (\ref{EPPFPKA}) for $\alpha=1/2$ and $\gamma=1$ one obtains
\begin{equation}
\label{EPPFIG}
p(n_1,\dots, n_k)= \frac{e^\delta \delta^k 2^{n-k}}{\Gamma(n)} \prod_{j=1}^k \left(\frac 12\right)_{n_j -1\uparrow} \int_0^\infty \lambda^{n-1} \frac{e^{-\delta(1+2\lambda)^{\frac 12}}}{\left( 1+ 2\lambda\right)^{n-\frac k2 }}d\lambda.
\end{equation}
With some manipulations, and having at hand the definition of incomplete Gamma function, i.e. $\Gamma(a;x)=\int_{x}^{\infty} t^{a-1} e^{-t}dt$, it is easy to see that (\ref{EPPFIG}) reduces to formula (A1) in Appendix A.4 of  Lijoi, Mena and Pr\"unster (2005), and results in Proposition 3 arise by specializing (\ref{weights}) and (\ref{vnk}) for $\alpha=1/2$.\\\\\\
It is worth to notice that the EPPF in (\ref{EPPFPKA}) defines an infinite {\it Gibbs partition of type $\alpha$},  (see Gnedin and Pitman, 2005), namely, for all $1 \leq k \leq n$, and all compositions $(n_1, \dots, n_k)$ of $n$, and for each $n\geq 1$, it has {\it Gibbs product form} i.e.
\begin{equation}
\label{GIBBSA}
p(n_1,\dots, n_k)= V_{n,k}\prod_{j=1}^k {W_{n_j}},
\end{equation}
for $W=(W_j)$ non-negative weights and
\begin{eqnarray}
&&V_{n,k} =\frac{e^{\delta\gamma} \delta^k \alpha^k 2^n }{\Gamma(n)} \int_{0}^{\infty} \lambda^{n-1} \frac{e^{-\delta (\gamma^{\frac 1\alpha}+2\lambda)^\alpha}}{(\gamma^{\frac 1\alpha}+2\lambda)^{n-k\alpha}}d\lambda,\nonumber
\end{eqnarray}
and it is of {\it type $\alpha$}, i.e. $W_{n_j}={(1-\alpha)_{n_j-1 \uparrow}}$. 
This is in line with a result stated in Pitman (2003) (cfr. Th. 8) and proved in Gnedin and Pitman (2005), according to which: {\it a)} an infinite exchangeable partitions $\Pi$ of $\mathbb{N}$ has EPPF in Gibbs form (\ref{GIBBSA}) if and only if $W_{n_j}={(1-\alpha)_{n_j-1 \uparrow}}$, for some $\alpha \in (-\infty ,1)$ and {\it b)} for $\alpha \in (0,1)$ this characterizes EPPFs induced by $PK(\rho_\alpha, \gamma)$ partition models, (cfr. Th. 12, item (iii) in Gnedin and Pitman, 2005).

\section{Distribution of the number of blocks}

In Bayesian nonparametric mixture modeling context it is usual to give to the distribution of the number of blocks in the partition induced by the prior,  the interpretation of a prior distribution for the number of components in the mixture model. Antoniak (1974) obtains the distribution of the number of component induced by a Dirichlet prior from the Ewens sampling formula, an equivalent of the EPPF for the Dirichlet case. From Gnedin and Pitman (2005) (cfr. eq. 10), an EPPF in Gibbs form induces the following distribution of the number of blocks $K_n$, by summation over all partitions  $\{A_1,\dots, A_k\}$ of $[n]$ with $k$ blocks, and $|A_j|=n_j$, $j=1,\dots,k$,  
\begin{equation}
\label{bloc}
Pr(K_n=k)=V_{n,k} B_{n,k}(W),
\end{equation}
where $B_{n,k}(W)= \sum_{\{A_1,\dots, A_k\}} \prod_{j=1}^k W_{|A_j|}$, is known as the partial Bell polynomial in the variables $W$.  A special form of (\ref{bloc}) for EPPFs in Gibbs form of type $\alpha$ can be easily derived to get the following result that don't need to be proved.\\\\
{\bf Proposition 2.} {\it A sample $(X_1,\dots,X_n)$ from a $PK(\rho_\alpha, \gamma_\alpha^e)$ model induces the following  distribution of the number of blocks $K_n$: 
\begin{equation}
\label{blocks}
Pr(K_n=k)=V_{n,k} S_\alpha (n,k)
\end{equation}
for 
\begin{eqnarray}
&&V_{n,k}=\frac{e^{\delta\gamma} \delta^k \alpha^k 2^n }{\Gamma(n)} \int_{0}^{\infty} \lambda^{n-1} \frac{e^{-\delta (\gamma^{\frac 1\alpha}+2\lambda)^\alpha}}{(\gamma^{\frac 1\alpha}+2\lambda)^{n-k\alpha}}d\lambda,\nonumber
\end{eqnarray}
and 
\begin{eqnarray}
&&S_{\alpha}(n, k):=B_{n,k}((1-\alpha)_{\bullet -1\uparrow}) = \frac{n!}{k!}\sum_{(n_1,\dots, n_k)} \prod_{j=1}^k \frac{1}{n_j!}(1-\alpha)_{n_j-1 \uparrow}\nonumber
\end{eqnarray}
where the sum extends over the space of all compositions $(n_1,\dots, n_k)$ of $n$. $S_\alpha(n,k)$ is known as the generalized Stirling number of the first kind, and it has the following explicit formula (cfr. Pitman, 2002)}.
\begin{eqnarray}
&&S_{\alpha}(n,k)=\frac{1}{\alpha^k k!}\sum_{j=1}^k (-1)^j {k \choose j} (-j\alpha)_{n \uparrow}\nonumber.
\end{eqnarray}\\\\
{\bf Example 3.}  The distribution of the number of components in a hierarchical Bayesian nonparametric mixture model under a N-IG prior obtained  in Lijoi, Mena and Pr\"unster (2005), Proposition 4.,  easily follows from (\ref{blocks}) for $\alpha =1/2$. In fact, from formula (127) in Pitman (2003), 
\begin{eqnarray}
&&B_{n,k}\left(\left(\frac 12 \right)_{\bullet -1 \uparrow} \right)=\frac{n!}{k!}\sum_{(n_1,\dots,n_k)}\prod_{j=1}^k
\frac{1}{n_j!}\left(\frac 12\right)_{n_j-1 \uparrow}= {2n-k-1 \choose n-1}\frac{\Gamma (n)}{\Gamma (k)}2^{2k-2n},\nonumber
\end{eqnarray}  
hence, for $\gamma=1$,
\begin{eqnarray}
&&Pr(K_n=k)=\frac{e^\delta \delta^k}{\Gamma(k)2^{n-k}}{2n-k-1 \choose n-1}  \int_0^\infty \lambda^{n-1} \frac{e^{-\delta(1+2\lambda)^{\frac 12}}}{\left( 1+ 2\lambda\right)^{n-\frac k2 }}d\lambda, \nonumber
\end{eqnarray}
which is easy to show that reduces to equation (9) in Lijoi, Mena and Pr\"unster (2005) by means of incomplete Gamma function substitutions.\\\\
In Section 6.1 of Pitman (2003) the concept of $\alpha${\it-diversity} has been introduced for a random partition $\Pi$ with ranked frequencies $(P_i)$ following a Poisson-Kingman model derived from an $\alpha$-stable law. \\\\
{\bf Definition 4.} [Pitman, 2003] An exchangeable partition $\Pi$ of the positive integers $\mathbb{N}$ has $\alpha$-diversity $S_\alpha$, if and only if there exists a random variable $S_\alpha$, with $0< S_\alpha <\infty$ a.s., such that, if $K_n$ is the number of blocks in the restriction of $\Pi$ to $[n]$, then 
\begin{equation}
\frac{K_n}{n^\alpha} \stackrel{a.s.}{\longrightarrow} S_\alpha {\mbox{\hspace{2cm} as $n \rightarrow \infty$}}.
\end{equation}
In Proposition 13, item $(i)$, Pitman  states that if $\Pi$ is a $PK(\rho_\alpha, \gamma)$ partition of $\mathbb{N}$ for some $\alpha \in (0,1)$, then $S_\alpha=T^{-\alpha}$, where $T$ has distribution $\gamma$. Hence, by an elementary transformation, the asymptotic distribution of ${K_n}/n^\alpha$ for a $PK(\rho_\alpha, \gamma_\alpha^e)$ partition, easily follows:\\\\
{\bf Proposition 3.} {\it Let $K_n$ be the number of blocks in a random  partition of $[n]$ induced by sampling from $P$, whose ranked atoms follow an exponentially tilted $PK(\rho_\alpha, \gamma_\alpha^e)$ model. Since $T$ has density 
\begin{eqnarray}
&&\gamma_\alpha^e(t)= exp\left\{\delta\gamma -\frac12 \gamma^{\frac 1\alpha} t\right\} f_{\alpha,\delta}(t)\nonumber,
\end{eqnarray}
for $f_{\alpha,\delta}(t)$ the density of an $\alpha$ stable law, then
\begin{eqnarray}
&&\frac{K_n}{n^\alpha} \stackrel{a.s.}{\longrightarrow} S_\alpha {\mbox{\hspace{2cm} as $n \rightarrow \infty$}}\nonumber
\end{eqnarray}
where $S_\alpha$ has density}
\begin{eqnarray}
\label{diversity}
f_{S_\alpha}(s)=exp\left\{{\delta\gamma -\frac12
\left(\frac{\gamma}{s}
\right)^{\frac 1\alpha}}\right\}\frac{f_{\alpha, \delta}(s^{-\frac 1\alpha})}{\alpha s^{\frac 1\alpha +1}}.
\end{eqnarray}\\\\
{\bf Example 4.} [Normalized inverse Gaussian process] For $\alpha=1/2$, $\gamma=1$ and exploiting the explicit form of the $1/2$-stable density, (\ref{diversity}) results:
\begin{eqnarray}
&&f_{S_{1/2}}(s)=\frac{\sqrt{2}\delta}{\sqrt{\pi}}\exp\left\{\delta -\frac 12 (\delta^2 s^2 + s^{-2})\right\}\nonumber.
\end{eqnarray}

\section*{References}
\newcommand{\bib}{\item \hskip-1.0cm}
\begin{list}{\ }{\setlength\leftmargin{1.0cm}}

\bib \textsc{Antoniak, C. E.} (1974) Mixtures of Dirichlet processes with applications to Bayesian nonparametric problems. {\it Ann. Statist.} 2, 1152-1174. 

\bib \textsc{Barndorff-Nielsen, O. E. and Shepard, N.} (2001)  Normal modified stable processes.  {\it Th. Probab. Math. Statist.}, 65, 1-19.

\bib \textsc{Brix, A.} (1999) Generalized Gamma measures and shot-noise Cox processes. {\it Adv. Appl. Probab.}, 31, 929--953.

\bib \textsc{Ferguson, T. S.} (1973)  A Bayesian analysis of some nonparametric problems. {\it Ann. Statist.}, 1, 209--230.

\bib \textsc{Gnedin, A. and Pitman, J. } (2005) {Exchangeable Gibbs partitions  and Stirling triangles} {\it arXiv:math. PR/0412494}

\bib \textsc{Hansen, B. and Pitman, J.} (2000) Prediction rules for exchangeable sequences related to species sampling. {\it Statistics \& Probability Letters}, 46, 251--256.

\bib \textsc{James, L. F.} (2002). Poisson process partition calculus with applications to exchangeable models and Bayesian Nonparametrics. {\it arXiv: math. ST/0205093}.

\bib \textsc{James, L.F., Lijoi, A. and Pr\"unster I.} (2005) Bayesian inference via classes of normalized random measures. {\it arXiv:math.ST/0503394}.

\bib \textsc{Kingman, J.F.C.} (1975) Random discrete distributions. {\it J. Roy. Statist. Soc. B}, 37, 1--22. 

\bib \textsc{Kingman, J.F.C} (1978) The representation of partition structure.  {\it J. London Math. Soc.} 2, 374--380.

\bib \textsc{Lijoi, A., Mena, R. and Pr\"unster, I.} (2005) Hierarchical mixture modeling with normalized Inverse-Gaussian priors. {\it JASA}, vol. 100, 1278--1291.

\bib \textsc{Pitman, J.} (1996) Some developments of the Blackwell-MacQueen urn scheme. In T.S. Ferguson, Shapley L.S., and MacQueen J.B., editors, {\it Statistics, Probability and Game Theory}, volume 30 of {\it IMS Lecture Notes-Monograph Series}, pages 245--267. Institute of Mathematical Statistics, Hayward, CA.

\bib \textsc{Pitman, J.} (2002) {\it Combinatorial Stochastic Processes}.  Lecture Notes for the Saint Flour summer school. Technical Report no.621 Dept. Statistics, U.C. Berkeley. To appear in Springer Lecture Notes in Mathematics. {\tt http://stat-www.berkeley.edu/users/pitman/621.ps.Z}

\bib \textsc{Pitman, J.} (2003) {Poisson-Kingman partitions}. In D.R. Goldstein, editor, {\it Science and Statistics: A Festschrift for Terry Speed}, volume 30 of Lecture Notes-Monograph Series, pages 1--34. Institute of Mathematical Statistics, Hayward, California.

\bib \textsc{Pitman, J. and Yor, M.} (1997) The two-parameter Poisson-Dirichlet distribution derived from a stable subordinator. {\it Ann. Probab.}, 25:855--900.

\bib \textsc{Seshadri, V. (1993)} {\it The inverse Gaussian distribution}. Oxford University Press, New York. 

\bib \textsc{Tweedie, M.} (1984). An index which distinguishes between some important exponential families. In J. Ghosh and J. Roy (Eds.), {\it Statistics: Applications and New Directions}: Proc. Indian Statistical Institute Golden Jubilee International Conference, pp. 579–-604.
\end{list}

\end{document}